\def\sgn{\operatorname{sgn}}
\def\Aut{\operatorname{Aut}}
\def\End{\operatorname{End}}
\def\Diff{\operatorname{Diff}}
\def\diff{\operatorname{diff}}
\def\Arg{\operatorname{Arg}}
\def\Diag{\operatorname{Diag}}
\def\Sp{\operatorname{Sp}}
\def\Tr{\operatorname{Tr}}
\def\SpQHS{\mathfrak{sp}_{{\mbox{\tiny{HS}}}}^{{\mbox{\tiny{Q}}}}}
\def\SpHS{\mathfrak{sp}_{\mbox{\tiny{HS}}}}

\documentclass[10pt,draft,reqno]{amsart}
\usepackage{amssymb}
     \makeatletter
     \def\section{\@startsection{section}{1}%
     \z@{.7\linespacing\@plus\linespacing}{.5\linespacing}%
     {\bfseries
     \centering
     }}
     \def\@secnumfont{\bfseries}
     \makeatother
\newtheorem{theorem}{Theorem}[section]
\newtheorem{lemma}[theorem]{Lemma}
\newtheorem{proposition}[theorem]{Proposition}
\newtheorem{corollary}[theorem]{Corollary}
\theoremstyle{definition}
\newtheorem{definition}[theorem]{Definition}
\theoremstyle{definition}
\newtheorem{notation}[theorem]{Notation}

\theoremstyle{remark}
\newtheorem{remark}[theorem]{Remark}
\numberwithin{equation}{section} 

\begin{document}

\title[$\Diff(S^{1})$ and $\Sp(\infty)$]
{Diffeomorphisms of the circle and Brownian motions on
an infinite-dimensional\\ symplectic group}

\author{Maria Gordina}
\thanks{* This research is partially supported by NSF
Grant DMS-0706784.}
\address{Department of Mathematics, University of Connecticut, Storrs, CT 06269, U.S.A. }
\email{gordina@math.uconn.edu}

\author{Mang Wu}
\thanks{* This research is partially supported by NSF
Grant DMS-0706784.}
\address{Department of Mathematics, University of Connecticut, Storrs, CT 06269, U.S.A. }
\email{mwu@math.uconn.edu}

\subjclass[2000] {Primary  60H07; Secondary  58J65, 60J65}

\keywords{$\Diff(S^{1})$, infinite-dimensional symplectic group,
Brownian motion}

\maketitle

\begin{abstract}
An embedding of the group $\Diff(S^{1})$ of orientation preserving
diffeomorphims of the unit circle $S^1$ into an infinite-dimensional
symplectic group, $\Sp(\infty)$, is studied. The authors prove that
this embedding is not surjective. A Brownian motion is constructed
on $\Sp(\infty)$. This study is motivated by recent work of
H.~Airault, S. Fang and P.~Malliavin.

\end{abstract}

\section{Introduction}
The group $\Diff(S^{1})$ of orientation preserving diffeomorphims of
the unit circle $S^1$ has been extensively studied for a long time.
One of the goals of the research has been to construct and study the
properties of a Brownian motion on this group. In \cite{AirMall2001}
H.~Airault and P.~Malliavin considered an embedding of $\Diff(S^1)$
into an infinite-dimensional symplectic group.

This group, $\Sp(\infty)$, can be represented as a certain
infinite-dimensional matrix group. For such matrix groups, the
method of\cite{Gordina2000a, Gordina2005a} can be used to construct
a Brownian motion living in the group. This construction relies on
the fact that these groups can be embedded into a larger Hilbert
space of Hilbert-Schmidt operators. We use the same method to
construct a Brownian motion on $\Sp(\infty)$. One of the advantages
of Hilbert-Schmidt groups is that one can associate an
infinite-dimensional Lie algebra to such a group, and this Lie
algebra is a Hilbert space. This is not the case with $\Diff(S^1)$,
as an infinite-dimensional Lie algebra associated with $\Diff(S^1)$
is not a Hilbert space with respect to the inner product compatible
with the symplectic structure on $\Diff(S^1)$.

In the current paper, we describe in detail the embedding of
$\Diff(S^1)$ into $\Sp(\infty)$, and construct a Brownian motion on
$\Sp(\infty)$. Our motivation comes from an attempt to use this
embedding to better understand Brownian motion in $\Diff(S^1)$ as
studied by H.~Airault, S. Fang and P.~Malliavin in a number of
papers (e.g. \cite{AirMall2001, AirMall2006, Fang2002,
FangLuo2007}). One of the main results of the paper is Theorem
\ref{t.4.6}, where we describe the embedding of $\Diff(S^1)$ into
$\Sp(\infty)$ and prove that the map is not surjective. Theorem
\ref{Main2} gives the construction of a Brownian motion on
$\Sp(\infty)$. In order for this Brownian motion to live in the
group we are forced to choose a non-$\operatorname{Ad}$-invariant
inner product on the Lie algebra of $\Sp(\infty)$. This fact has a
potential implication for this Brownian motion not to be
quasi-invariant for the appropriate choice of the Cameron-Martin
subgroup of $\Sp(\infty)$. This is in contrast to results in
\cite{AirMall2006}. The latter can be explained by the fact that the
Brownian motion we construct in Section \ref{s.6} lives in a
subgroup of $\Sp(\infty)$ whose Lie algebra is much smaller than the
full Lie algebra of $\Sp(\infty)$.

\section{The spaces $H$ and $\mathbb{H}_\omega$}

\begin{definition}\label{def.Omega}
Let $H$ be the space of complex-valued $C^\infty$ functions on the
unit circle $S^1$ with the mean value $0$. Define a bilinear form
$\omega$ on $H$ by
\[
\omega(u,v)=\frac{1}{2\pi}\int_0^{2\pi}uv'd\theta,
\hspace{.2in}\text{ for any } u,v\in H.
\]
\end{definition}

\begin{remark} By using integration by parts, we see that the form
$\omega$ is anti-symmetric, that is,  $\omega(u,v)=-\omega(v,u)$ for
any $u, v \in H$.
\end{remark}
Next we define an inner product $(\cdot,\cdot)_\omega$ on $H$ which
is compatible with the form $\omega$. First, we introduce a complex
structure on $H$, that is, a linear map $J$ on $H$ such that
$J^2=-id$. Then the inner product is defined by $(u,v)_\omega=\pm
\omega(u,J\bar{v})$, where the sign depends on the choice of $J$.
The complex structure $J$ in this context is called the Hilbert
transform.

\begin{definition}\label{def.L2Space}
Let $\mathbb{H}_0$ be the Hilbert space of complex-valued $L^2$
functions on $S^1$ with the mean value $0$ equipped with the inner
product
\[
(u,v) = \frac{1}{2\pi}\int_0^{2\pi}u\bar{v}d\theta, \hspace{.2in}
\text{ for any }  u,v\in\mathbb{H}_0.
\]
\end{definition}


\begin{notation}\label{Decomp}Denote $\hat{e}_n=e^{in\theta}, n \in \mathbb{Z}\backslash\{0\}$, and
$ \mathcal{B}_H =\left\{\hat{e}_n, \ n \in \mathbb{Z}\backslash\{0\}
\right\}.$ Let $\mathbb{H}^+$ and $\mathbb{H}^-$ be the closed
subspaces of $\mathbb{H}_0$ spanned by $\{\hat{e}_n:n>0\}$ and
$\{\hat{e}_n:n<0\}$, respectively. By $\pi^+$ and $\pi^-$ we denote
the projections of $\mathbb{H}_0$ onto subspaces $\mathbb{H}^+$ and
$\mathbb{H}^-$, respectively. For $u\in\mathbb{H}_0$, we can write
$u=u_++u_-$, where $u_+=\pi^+(u)$ and $u_-=\pi^-(u)$.
\end{notation}

\begin{definition}\label{def.J}
Define the \textbf{Hilbert transformation} $J$ on $\mathcal{B}_H$ by
\[
J:\hat{e}_n \mapsto i\sgn(n)\hat{e}_n
\]
where $\sgn(n)$ is the sign of $n$, and then extended by linearity
to $\mathbb{H}_0$.
\end{definition}

\begin{remark}
In the above definition, $J$ is defined on the space
$\mathbb{H}_{0}$. We need to address the issue whether it is
well--defined on the \emph{subspace} $H$. That is, if $J(H)\subseteq
H$. We will see that if we modify the space $H$ a little bit, for
example, if we let $C_0^1(S^1)$ be the space of complex-valued $C^1$
functions on the circle with mean value zero, then $J$ is \emph{not}
well--defined on $C_0^1(S^1)$. This problem really lies in the heart
of Fourier analysis. To see this, we need to characterize $J$ by
using the Fourier transform.
\end{remark}

\begin{notation}
For $u\in\mathbb{H}_0$, let $\mathcal{F}:u\mapsto\hat{u}$ be the
\textbf{Fourier transformation} with $\hat{u}(n)=(u,\hat{e}_n)$. Let
$\hat{J}$ be a transformation on $l^2(\mathbb{Z}\backslash\{0\})$
defined by $\big(\hat{J}\hat{u}\big)(n)=i\sgn(n)\hat{u}(n)$ for any
$\hat{u}\in l^2(\mathbb{Z}\backslash\{0\})$.
\end{notation}

The Fourier transformation $\mathcal{F}:\mathbb{H}_0 \to
l^2(\mathbb{Z}\backslash\{0\})$ is an isomorphism of Hilbert spaces,
and  $J=\mathcal{F}^{-1} \circ \hat{J} \circ \mathcal{F}$.

\begin{proposition}\label{JWell}
The Hilbert transformation $J$ is well--defined on $H$, that is
$J(H)\subseteq H$.
\end{proposition}

\begin{proof}
The key of the proof is the fact that functions in  $H$ can be
completely characterized by their Fourier coefficients. To be
precise, let $u\in\mathbb{H}_0$ be continuous. Then $u$ is in
$C^\infty$ if and only if $\lim_{n\to\infty}n^k\hat{u}(n)=0$ for any
$k\in\mathbb{N}$. From this fact, it follows immediately that $J$ is
well--defined on $H$, because $J$ only changes the signs of the
Fourier coefficients of a function $u\in H$.

For completeness of exposition, we give a proof of this
characterization. Though this is probably a standard fact in Fourier
analysis, we found a proof (in \cite{KatznelsonBook}) of only one
direction.

We first assume that $u$ is $C^\infty$. Then
$u(\theta)=u(0)+\int_0^\theta u'(t)dt$. So
\begin{align*}
\hat{u}(n)
&=\frac{1}{2\pi}\Big(\int_0^{2\pi}\int_0^{2\pi}u'(t)\chi_{[0,\theta]}
dt\Big) e^{-in\theta}d\theta
=\frac{1}{2\pi}\int_0^{2\pi}\Big(\int_t^{2\pi}e^{-in\theta}d\theta\Big)u'(t)dt\\
&=-\frac{1}{2\pi in}\int_0^{2\pi}u'(t)-u'(t)e^{-int}dt
=\frac{\widehat{u'}(n)}{in},
\end{align*}
where we have used Fubini's theorem and the continuity of $u'$. Now,
$u'$ is itself $C^\infty$, so we can apply the procedure again. By
induction, we get
$\hat{u}(n)=\frac{\widehat{u^{(k)}}(n)}{\left(in\right)^{k}}$. But
from the general theory of Fourier analysis,
$\widehat{u^{(k)}}(n)\to 0$ as $n\to\infty$. Therefore
$n^k\hat{u}(n)\to 0$ as $n\to\infty$.

Conversely, assume $u$ is such that for any $k$, $n^k\hat{u}(n)\to
0$ as $n\to\infty$. Then the Fourier series of $u$ converges
uniformly. Also by assumption that $u$ is continuous, the Fourier
series converges to $u$ for all $\theta\in S^1$ (see Corollary I.3.1
in \cite{KatznelsonBook}). So we can write $u(\theta) = \sum_{n\neq
0} \hat{u}(n) e^{in\theta}$. Fix a point $\theta\in S^1$, then
\[
u'(\theta)
= \left.\frac{d}{dt}\right|_{t=\theta}
\sum_{n\neq 0} \hat{u}(n) e^{int}
=\lim_{t\to\theta}\lim_{N\to\infty}
\sum_{n=-N}^{N} \hat{u}(n)
\frac{e^{int}-e^{in\theta}}{t-\theta}.
\]
Note that the derivatives of $\cos nt$ and $\sin nt$ are all bounded
by $|n|$. So by the mean value theorem, $|\cos nt-\cos n\theta|\le
|n||t-\theta|$, and $|\sin nt-\sin n\theta|\le |n||t-\theta|$. So
\[
\Big|\frac{e^{int}-e^{in\theta}}{t-\theta}\Big| \le 2|n|,
\hspace{.2in}\text{ for any } t,\theta\in S^1.
\]
Therefore, by the growth condition on the Fourier coefficients
$\hat{u}$, we have
\[
\lim_{N\to\infty}\sum_{n=-N}^{N}
\hat{u}(n) \frac{e^{int}-e^{in\theta}}{t-\theta}
\]
converges at the fixed $\theta\in S^1$ and the convergence is uniform in $t\in S^1$.
Therefore we can interchange the two limits, and obtain
\[
\Big(\sum_{n\neq 0} \hat{u}(n) e^{in\theta}\Big)'
=\sum_{n\neq 0} \hat{u}(n) in e^{in\theta},
\]
which means we can differentiate term by term. So the Fourier
coefficients of $u'$ are given by $\hat{u'}(n)=in\hat{u}(n)$.
Clearly, $\hat{u'}$ satisfies the same condition as $\hat{u}$:
$n^k\hat{u'}(n)\to 0$ as $n\to\infty$. By induction, $u$ is
$j$-times differentiable for any $j$. Therefore, $u$ is in
$C^\infty$.
\end{proof}

\begin{proposition}
Let $C_0^1(S^1)$ be the space of complex-valued $C^1$ functions on
the circle with the mean value zero. Then the Hilbert transformation
$J$ is \emph{not} well--defined on $C_0^1(S^1)$, i.e.,
$J(C_0^1(S^1))\nsubseteq C_0^1(S^1)$.
\end{proposition}

\begin{proof}
Let $C(S^1)$ be the space of continuous functions on the circle. In
 \cite{KatznelsonBook}, it is shown that there
exists a function in $C(S^1)$ such that the corresponding Fourier
series does not converges \emph{uniformly} \cite[Theorem
II.1.3]{KatznelsonBook}, and therefore there exists an $f\in C(S^1)$
such that $Jf\notin C(S^1)$ \cite[Theorem II.1.4]{KatznelsonBook}.
Now take $u=f-f_0$ where $f_0$ is the mean value of $f$. Then $u$ is
a continuous function on the circle with the mean value zero, and
$Ju$ is \emph{not} continuous.

Using Notation \ref{Decomp} let us write $u=u_+ + u_-$. Then we can
use the relation
\[
iu+Ju=2iu_+ \hspace{.1in}\mbox{and}\hspace{.1in} iu-Ju=2iu_-.
\]
to see that $u_+$ and $u_-$ are \emph{not} continuous. Integrating
$u=u_+ + u_-$, we have
\[
\int_0^t u(\theta)d\theta
=\int_0^t u_+(\theta)d\theta
+\int_0^t u_-(\theta)d\theta.
\]
Denote the three functions in the above equation by $v,v_1,v_2$.
By theorem I.1.6 in \cite{KatznelsonBook},
\[
\hat{v}(n)=\frac{\hat{u}(n)}{in},
\hspace{.1in}\mbox{and}\hspace{.1in}
\hat{v_1}(n)=\frac{\hat{u}_+(n)}{in},
\hat{v_2}(n)=\frac{1}{in}\hat{u}_-(n) \mbox{ for } n\neq 0.
\]

Let $g=v-v_0$ where $v_0$ is the mean value of $v$. Then $g\in C_0^1(S^1)$.
Write $g=g_+ + g_-$ \ref{Decomp}.
Then $g_+=v_1-(v_1)_0$ and $g_-=v_2-(v_2)_0$ where
$(v_1)_0$ and $(v_2)_0$ are the mean values of $v_1$ and $v_2$ respectively.
Then $g_+,g_-\notin C_0^1(S^1)$ since $v_1'=u_+,v_2'=u_-$ are \emph{not} continuous.

By the relation
\[
ig+Jg=2ig_+ \hspace{.1in}\mbox{and}\hspace{.1in} ig-Jg=2ig_-,
\]
we see that $Jg\notin C_0^1(S^1)$.
\end{proof}

\begin{notation}\label{nota.Inner}
Define an $\mathbb{R}$-bilinear form $(\cdot,\cdot)_\omega$ on $H$
by
\[
(u,v)_\omega=-\omega(u,J\bar{v}) \hspace{.2in}\text{ for any }
u,v\in H.
\]

\end{notation}

\begin{proposition}
$(\cdot,\cdot)_\omega$ is an inner product on $H$.
\end{proposition}

\begin{proof}
We need to check that $(\cdot,\cdot)_\omega$ satisfies the following
properties (1) $(\lambda u,v)_\omega=\lambda(u,v)_\omega$ for
$\lambda\in\mathbb{C}$; (2) $(v,u)_\omega=\overline{(u,v)_\omega}$;
(3) $(u,u)_\omega>0$ unless $u=0$.

(1) for $\lambda\in\mathbb{C}$,
\[
(\lambda u,v)_\omega =-\omega(\lambda u,J\bar{v})
=-\lambda\cdot\omega(u,J\bar{v}) =\lambda\cdot(u,v)_\omega.
\]

To prove (2) and (3), we need some simple facts:
$H^+=\pi^+(H)\subseteq H$ and $H^-=\pi^-(H)\subseteq H$, and
$H=H^+\oplus H^-$. If $u\in H^+,v\in H^-$, then $(u,v)=0$. If $u\in
H^+$, then $\bar{u}\in H^-, Ju=iu, Ju\in H^+$. If $u\in H^-$, then
$\bar{u}\in H^+, Ju=-iu, Ju\in H^-$. $J\bar{u}=\overline{Ju}$.
$\widehat{u'}(n)=in\hat{u}(n)$. In particular, if $u\in H^+$, then
$u'\in H^+$; if $u\in H^-$, then $u'\in H^-$.

(2)
By definition,
\begin{align*}
& (v,u)_\omega=-\omega(v,J\bar{u})=\omega(J\bar{u},v)
=\frac{1}{2\pi}\int (J\bar{u})v' d\theta \notag
\\ &
\overline{(u,v)_\omega}=-\overline{\omega(u,J\bar{v})}
=\overline{\omega(J\bar{v},u)} =\frac{1}{2\pi}\int
\overline{J\bar{v}}\bar{u}' d\theta =\frac{1}{2\pi}\int (Jv)\bar{u}'
d\theta. \label{e.2.1}
\end{align*}
Write $u=u_+ + u_-$ and $v=v_+ + v_-$ as in Notation \ref{Decomp}.
Using the above fact, we can show that the above two quantities are
equal to each other.

(3) Write $u=u_+ + u_-$, then
\[
(u,u)_\omega
=\frac{1}{2\pi}\int (-i\overline{u_+}u_+' + i\overline{u_-}u_-') d\theta
=\sum_{n\neq 0}|n||\hat{u}(n)|^2.
\]
Therefore, $(u,u)_\omega>0$ unless $u=0$.
\end{proof}

\begin{definition}\label{def.OmegaBasis}
Let $\mathbb{H}_\omega$ be the completion of $H$ under the norm
$\|\cdot\|_\omega$ induced by the inner product
$(\cdot,\cdot)_\omega$. Define
\[
\mathcal{B}_\omega=
\left\{\tilde{e}_n=\frac{1}{\sqrt{n}}e^{in\theta},
n>0 \right\} \cup \left\{\tilde{e}_n=\frac{1}{i\sqrt{|n|}}e^{in\theta},
n<0 \right\}.
\]
\end{definition}

\begin{remark}\label{Completion}
$\mathbb{H}_\omega$ is a Hilbert space. Also the norm
$\|\cdot\|_\omega$ induced by the inner product
$(\cdot,\cdot)_\omega$ is \emph{strictly} stronger than the norm
$\|\cdot\|$ induced by the inner product $(\cdot,\cdot)$. So
$\mathbb{H}_\omega$ can be identified as a \emph{proper} subspace of
$\mathbb{H}_0$. The inner product $(\cdot,\cdot)_\omega$ or the norm
induced by it is sometimes called the $H^{1/2}$ metric or the
$H^{1/2}$ norm on the space $H$.

One can verify that $\mathcal{B}_\omega$ is an orthonormal basis of
$\mathbb{H}_\omega$. From the definition of the inner product
$(\cdot,\cdot)_\omega$, we have the relation
$\omega(u,v)=(u,\overline{Jv})_\omega$ for any $u,v\in H$. This can
be used to \emph{extend} the form $\omega$ to $\mathbb{H}_\omega$.

Finally, from the non--degeneracy of the inner
product $(\cdot,\cdot)_\omega$, we see that the form
$\omega(\cdot,\cdot)$ on $\mathbb{H}_\omega$ is also
non--degenerate.
\end{remark}

\section{An infinite-dimensional symplectic group}

\begin{definition}\label{def.BoundedOp}
Let $B(\mathbb{H}_\omega)$ be the space of \textbf{bounded
operators} on $\mathbb{H}_\omega$ equipped with the operator norm.
For an operator $A\in B(\mathbb{H}_\omega)$

\begin{enumerate}

\item
suppose $\bar{A}$ is an operator on $\mathbb{H}_\omega$ satisfying
$\bar{A}u=\overline{A\bar{u}}$ for any $u\in \mathbb{H}_\omega$,
then $\bar{A}$ is the \textbf{conjugate} of $A$;

\item suppose $A^\dag$ is an operator on $\mathbb{H}_\omega$ satisfying
$(Au,v)_\omega=(u,A^\dag v)_\omega$ for any $u,v\in\mathbb{H}_\omega$, then $A^\dag$ is the
\textbf{adjoint} of $A$;

\item  then $A^T=\bar{A}^\dag$ is the \textbf{transpose} of $A$;

\item  suppose $A^\#$ is an operator on $\mathbb{H}_\omega$ satisfying $\omega(Au,v)=\omega(u,A^\# v)$ for any
$u,v\in\mathbb{H}_\omega$, then $A^\#$ is the \textbf{symplectic
adjoint} of $A$.

\item $A$ is said to \textbf{preserve the form} $\omega$ if
$\omega(Au,Av)=\omega(u,v)$ for any $u,v\in\mathbb{H}_\omega$.
\end{enumerate}
\end{definition}

In the orthonormal basis $\mathcal{B}_\omega$, an operator $A\in
B(\mathbb{H}_\omega)$ can be represented by an infinite-dimensional
matrix, still denoted by $A$, with $(m,n)$th entry equal to
$A_{m,n}=(A\tilde{e}_n,\tilde{e}_m)_\omega$.

\begin{remark}\label{NegativeIdx}
If we represent an operator $A\in B(\mathbb{H}_\omega)$ by a matrix
$\{A_{m,n}\}_{m,n\in\mathbb{Z}\backslash\{0\}}$,
the indices $m$ and $n$ are allowed to be both positive and negative following
Definition \ref{def.OmegaBasis} of $\mathcal{B}_\omega$.
\end{remark}

The next proposition collects some simple facts about operations on
$B(\mathbb{H}_\omega)$ introduced in Definition \ref{def.BoundedOp}.

\begin{proposition}\label{SomeFacts1}
Let $A,B\in B(\mathbb{H}_\omega)$. Then
\begin{enumerate}
\item
$\overline{\tilde{e}_n}=i\tilde{e}_{-n}$,
$J\tilde{e}_n=i\sgn(n)\tilde{e}_n$,
$(\tilde{e}_n)'=in\tilde{e}_n$;

\item
$(\bar{A})_{m,n}=\overline{A_{-m,-n}}$;

\item
$(A^\dag)_{m,n}=\overline{A_{n,m}}$;

\item
$\bar{A}^\dag=\overline{A^\dag}$,
and
$(A^T)_{m,n}=A_{-n,-m}$;

\item
if $A=\bar{A}$, then $(A^\#)_{m,n}=\sgn(mn)\overline{A_{n,m}}$;

\item
$\overline{AB}=\bar{A}\bar{B}$,
$(AB)^\dag=B^\dag A^\dag$,
$(AB)^T=B^T A^T$,
$(AB)^\#=B^\# A^\#$;

\item
If $A$ is invertible,
then $\bar{A},A^T,A^\dag,A^\#$ are all invertible,
and $(\bar{A})^{-1}=\overline{A^{-1}}$,
$(A^T)^{-1}=(A^{-1})^T$,
$(A^\dag)^{-1}=(A^{-1})^\dag$,
$(A^\#)^{-1}=(A^{-1})^\#$;

\item
$(\pi^+)_{m,n}=\frac{1}{2}(\delta_{mn}+\sgn(m)\delta_{mn})$,
$(\pi^-)_{m,n}=\frac{1}{2}(\delta_{mn}-\sgn(m)\delta_{mn})$,
$\overline{\pi^+}=\pi^-$, $\overline{\pi^-}=\pi^+$,
$(\pi^+)^T=\pi^-$, $(\pi^-)^T=\pi^+$,
$(\pi^+)^\dag=\pi^+$, $(\pi^-)^\dag=\pi^-$;

\item
$J_{m,n}=i\sgn(m)\delta_{mn}$, $\bar{J}=J$, $J=i(\pi^+-\pi^-)$,
$J^T=-J$, $J^\dag=-J$, $J^2=-id$;

\item
$(A^\#)_{m,n}=\sgn(mn)A_{-n,-m}$.
\end{enumerate}
\end{proposition}

\begin{proof}
All of these properties can be checked by straight forward
calculations. We only prove (10).
\begin{align*}
&
(A^\#)_{m,n}
=(A^\#\tilde{e}_n,\tilde{e}_m)_\omega
=-\omega(A^\#\tilde{e}_n,J\overline{\tilde{e}_m})
=\omega(J\overline{\tilde{e}_m},A^\#\tilde{e}_n)
\\&
=\omega(AJ\overline{\tilde{e}_m},\tilde{e}_n)
=-\omega(\tilde{e}_n,AJ\overline{\tilde{e}_m})
=-\omega(\tilde{e}_n,J(-J)AJ\overline{\tilde{e}_m})
\\&
=-\omega(\tilde{e}_n,J\overline{(-J\bar{A}J\tilde{e}_m)}),
\end{align*}
where in the last equality we used property (6),
$\overline{AB}=\bar{A}\bar{B}$, and  property (9), $\bar{J}=J$, so
that
$\overline{-J\bar{A}J\tilde{e}_m}=-\bar{J}\bar{\bar{A}}\bar{J}\overline{\tilde{e}_m}
=-JAJ\overline{\tilde{e}_m}$. Therefore,
\begin{align*}
&
(A^\#)_{m,n}
=-\omega(\tilde{e}_n,J\overline{(-J\bar{A}J\tilde{e}_m)})
=(\tilde{e}_n,-J\bar{A}J\tilde{e}_m)_\omega
=-(\tilde{e}_n,J\bar{A}J\tilde{e}_m)_\omega
\\&
=-(J^\dag\tilde{e}_n,\bar{A}J\tilde{e}_m)_\omega
=-(-J\tilde{e}_n,\bar{A}J\tilde{e}_m)_\omega
=(i\sgn(n)\tilde{e}_n,\bar{A}i\sgn(m)\tilde{e}_m)_\omega
\\&
=\sgn(mn)(\tilde{e}_n,\bar{A}\tilde{e}_m)_\omega
=\sgn(mn)\overline{(\bar{A}\tilde{e}_m,\tilde{e}_n)_\omega}
=\sgn(mn)\overline{(\bar{A})_{n,m}}
\\&
=\sgn(mn)A_{-n,-m}.
\end{align*}
\end{proof}

\begin{notation}\label{nota.BlockMatrix}
For $A\in B(\mathbb{H}_\omega)$, let $a=\pi^+ A \pi^+$, $b=\pi^+ A
\pi^-$, $c=\pi^- A \pi^+$, and $d=\pi^- A \pi^-$, where
$a:\mathbb{H}_\omega^+\to\mathbb{H}_\omega^+$,
$b:\mathbb{H}_\omega^-\to\mathbb{H}_\omega^+$,
$c:\mathbb{H}_\omega^+\to\mathbb{H}_\omega^-$,
$d:\mathbb{H}_\omega^-\to\mathbb{H}_\omega^-$. Then $A=a+b+c+d$ can
be represented as the following block matrix
\[
\left(
\begin{array}{ll}
a &b\\
c &d
\end{array}
\right).
\]
\end{notation}

If $A, B\in B(\mathbb{H}_\omega)$, then the block matrix
representation for $AB$ is exactly the multiplication of block
matrices for $A$ and $B$.

\begin{proposition}
Suppose $A\in B(\mathbb{H}_\omega)$ with the matrix $\{A_{m,n}\}_{m,
n \in \mathbb{Z}\backslash\{0\}}$. Then the following are equivalent
\begin{enumerate}
\item
$A=\bar{A}$;

\item
if $u=\bar{u}$, then $Au=\overline{Au}$;

\item
$A_{m,n}=\overline{A_{-m,-n}}$ (\ref{NegativeIdx});

\item
as a block matrix, $A$ has the form
$
\left(
\begin{array}{ll}
a & b\\
\bar{b} &\bar{a}
\end{array}
\right).
$
\end{enumerate}
\end{proposition}

\begin{proof}
Equivalence of (1), (3) and (4) follows from
Proposition\ref{SomeFacts1} and Notation\ref{nota.BlockMatrix}.
First we  show that (1) is equivalent to (2).

[(1)$\Longrightarrow$(2)].
If $u=\bar{u}$, then $Au=\bar{A}u=\overline{A\bar{u}}=\overline{Au}$.

[(2)$\Longrightarrow$(1)]. Let
$u=\tilde{e}_n+\overline{\tilde{e}_n}$, and
$v=\tilde{e}_{-n}+\overline{\tilde{e}_{-n}}$. Then $u,v$ are
real-valued functions on the circle. Using Proposition
\ref{SomeFacts1} we have $\overline{\tilde{e}_n}=i\tilde{e}_{-n}$,
and therefore $Au=\overline{Au}$ and $Av=\overline{Av}$ imply
\begin{align*}
& A\tilde{e}_n + iA\tilde{e}_{-n} = \overline{A\tilde{e}_n} -
i\overline{A\tilde{e}_{-n}}
\\
& A\tilde{e}_n - iA\tilde{e}_{-n} = -\overline{A\tilde{e}_n} -
i\overline{A\tilde{e}_{-n}}.
\end{align*} Solving the above two
equations for $A\tilde{e}_n$, we have
\[
A\tilde{e}_n=-i\overline{A\tilde{e}_{-n}}
=\overline{A\overline{\tilde{e}_n}}=\bar{A}\tilde{e}_n
\]
with this being true for any $n\neq 0$, and so $A=\bar{A}$.
\end{proof}

\begin{proposition}\label{prop.Preserve}
Let $A\in B(\mathbb{H}_\omega)$.
The following are equivalent:

\begin{enumerate}
\item
$A$ preserves the form $\omega$;

\item
$\omega(Au,Av)=\omega(u,v)$ for any $u,v\in\mathbb{H}_\omega$;

\item
$\omega(A\tilde{e}_m,A\tilde{e}_n)=\omega(\tilde{e}_m,\tilde{e}_n)$
for any $m,n\neq 0$;

\item
$A^TJA=J$;

\item
$\sum_{k\neq 0}\sgn(mk)A_{k,m}A_{-k,-n}=\delta_{m,n}$ for any $m,n\neq0$.
\end{enumerate}
If we further assume that $A=\bar{A}$, then the following two are
equivalent to the above:

\begin{enumerate}
\item[(I)]
$a^T\bar{a}-b^\dag b=\pi^-$ and $a^T\bar{b}-b^\dag a=0$;

\item[(II)]
$\sum_{k\neq 0}\sgn(mk)A_{k,m}\overline{A_{k,n}}=\delta_{m,n}$ for any $m,n\neq0$.
\end{enumerate}

\end{proposition}

\begin{proof}
Equivalence of (1),(2) and (3) follows directly from Definition
\ref{def.BoundedOp}. Let us check the equivalency of (2) and (4).
First assume that (2) holds. By Remark \ref{Completion} we have
$\omega(u,v)=(u,J\bar{v})_\omega$ , and therefore
\[
\omega(Au,Av)=(Au,J\overline{Av})_\omega=(u,A^\dag J\overline{Av})_\omega.
\]
By assumption, $\omega(Au,Av)=\omega(u,v)$ for any
$u,v\in\mathbb{H}_\omega$. So by the non-degeneracy of the inner
product $(\cdot,\cdot)_\omega$, we have $A^\dag
J\overline{Av}=J\bar{v}$ for any $v\in\mathbb{H}_\omega$. By
definition of $\bar{A}$, we have $\overline{Av}=\bar{A}\bar{v}$. So
$A^\dag J\bar{A}\bar{v}=J\bar{v}$ for any $v\in\mathbb{H}_\omega$,
or $A^\dag J\bar{A}=J$. Taking conjugation of both sides and using
$\bar{J}=J$, we see that $A^T JA=J$.

Every step above is reversible, therefore we have implication in the
other direction as well.

Now we check the equivalency of (3) and (5). First, by Remark
\ref{Completion} $\omega(u,v)=(u,J\bar{v})_\omega$
 and Proposition \ref{SomeFacts1}
\[
\omega(\tilde{e}_m,\tilde{e}_n)
=(\tilde{e}_m,J\overline{\tilde{e}_n})_\omega
=-\sgn(m)\delta_{m,-n}.
\]
On the other hand, by the continuity of the form
$\omega(\cdot,\cdot)$ in both variables, we have
\begin{align*}
&
\omega(A\tilde{e}_m,A\tilde{e}_n)
=\omega\Big(\sum_k A_{k,m}\tilde{e}_k,\sum_k A_{l,n}\tilde{e}_l\Big)
\\ &
=\sum_{k,l}A_{k,m}A_{l,n}(-\sgn(k))\delta_{k,-l} =-\sum_k
\sgn(k)A_{k,m}A_{-k,n}.
\end{align*}
Now assuming
$\omega(A\tilde{e}_m,A\tilde{e}_n)=\omega(\tilde{e}_m,\tilde{e}_n)$,
we have
\[
-\sum_k \sgn(k)A_{k,m}A_{-k,n}=-\sgn(m)\delta_{m,-n},
\text{ for any } m,n\neq 0.
\]
By multiplying by $\sgn(m)$ both sides, and replacing $-n$ with $n$,
we get (5). Conversely, note that every step above is reversible,
therefore we have implication in the other direction.

We have proved equivalence of (1)-(5). Now assume $A=\bar{A}$. To
prove equivalence of (4) and (I), just notice that as block
matrices, $A,A^T$ and $J$ have the form
\[
\left(
\begin{array}{ll}
a & b\\
\bar{b} & \bar{a}
\end{array}
\right),
\hspace{.2in}
\left(
\begin{array}{ll}
a^\dag & b^T\\
b^\dag & a^T
\end{array}
\right),
\hspace{.1in}\mbox{and}\hspace{.1in}
i\left(
\begin{array}{ll}
\pi^+ & 0\\
0 & -\pi^-
\end{array}
\right).
\]
Equivalence of (5) and (II) follows 
from the relation
$A_{-k,-n}=\overline{A_{k,n}}$.
\end{proof}

\begin{proposition}\label{prop.Invertible}
Let $A\in B(\mathbb{H}_\omega)$.
If $A$ preserves the form $\omega$, then the following are equivalent:
\begin{enumerate}
\item
$A$ is invertible.

\item
$AJA^T=J$.

\item
$A^T$ preserves the form $\omega$.

\item
$\sum_k \sgn(mk)A_{m,k}A_{-n,-k}=\delta_{m,n}$ for any $m,n\neq0$.
\end{enumerate}

If we further assume that $A=\bar{A}$, then the following are
equivalent to the above:

\begin{enumerate}
\item[(I)]
$\bar{a}a^T-\bar{b}b^T=\pi^-$ and $\bar{b}a^\dag-\bar{a}b^\dag=0$.

\item[(II)]
$\sum_k \sgn(mk)A_{m,k}\overline{A_{n,k}}=\delta_{m,n}$ for any $m,n\neq0$.
\end{enumerate}

\end{proposition}

\begin{proof} We will use several times the fact that if $A$ preserves $\omega$, then
$A^TJA=J$.

[(1)$\Rightarrow$(2)] Multiplying on the left by $(A^T)^{-1}$ and
multiplying on the right by $A^{-1}$ both sides, we get
$J=(A^T)^{-1}JA^{-1}$, and so $(A^{-1})^TJA^{-1}=J$. Taking inverse
of both sides, and using $J^{-1}=-J$, we have $A^TJA=J$.

[(2)$\Rightarrow$(1)]  As $J$ is injective, so is $A^TJA$, and
therefore $A$ is injective. On the other hand, by assumption
$AJA^T=J$. As $J$ is surjective, so $AJA^T$ is surjective too. This
implies that $A$ is surjective, and therefore $A$ is invertible.

Equivalence of (2) and (3) follows from  $(A^T)^T=A$ and Proposition
\ref{prop.Preserve}. Equivalence of (3) and (4) follows directly
from Proposition \ref{prop.Preserve} and the fact that
$(A^T)_{m,n}=A_{-n,-m}$.

Now assume that $A=\bar{A}$. Then equivalence of (3) and (I)can be
checked by using multiplication of block matrices as in the proof of
Proposition \ref{prop.Preserve}. Finally (4) is equivalent to (II)
as if $A=\bar{A}$, then $A_{-m,-n}=\overline{A_{m,n}}$.
\end{proof}

\begin{corollary}\label{AASharp}
Let $A\in B(\mathbb{H}_\omega)$ and $A=\bar{A}$.
Then the following are equivalent:
\begin{enumerate}
\item
$A$ preserves the form $\omega$ and is invertible;
\item
$A^\#A=A^\#A=id$;
\end{enumerate}
\end{corollary}

\begin{proof}
By Proposition \ref{SomeFacts1}
\begin{align*}
&(A^\#A)_{m,n}=\sum_{k\neq0}(A^\#)_{m,k}A_{k,n}
=\sum_{k\neq0}\sgn(mk)A_{k,n}\overline{A_{k,m}},\\
&(AA^\#)_{m,n}=\sum_{k\neq0}A_{m,k}(A^\#)_{k,n}
=\sum_{k\neq0}\sgn(nk)A_{m,k}\overline{A_{n,k}}.
\end{align*}
Therefore, by (II) in Proposition \ref{prop.Preserve} and (II) in
Proposition \ref{prop.Invertible} we have equivalence.
\end{proof}

\begin{definition}\label{def.2Norm}
Define a
(semi)norm $\|\cdot\|_2$ on $B(\mathbb{H}_\omega)$ such that for
$A\in B(\mathbb{H}_\omega)$, $\|A\|_2^2=\Tr(b^\dag b)=\Vert b
\Vert_{HS}$, where $b=\pi^+ A \pi^-$. That is, the norm $\|A\|_2$ is
just the Hilbert-Schmidt norm of the block $b$.
\end{definition}

\begin{definition}\label{def.SPInf}
An \textbf{infinite-dimensional symplectic group} $\Sp(\infty)$ is
the set of bounded operators $A$ on $H$ such that
\begin{enumerate}
\item $A$ is invertible;
\item $A=\bar{A}$;
\item $A$ preserves the form $\omega$;
\item $\|A\|_2<\infty$.
\end{enumerate}
\end{definition}

\begin{remark}\label{def.SPInf1}
If $A$ is a bounded operator on $H$, then $A$ can be extended to a
bounded operator on $\mathbb{H}_\omega$. Therefore, we can
equivalently define $\Sp(\infty)$ to be the set of operators $A\in
B(\mathbb{H}_\omega)$ such that
\begin{enumerate}
\item $A$ is invertible;
\item $A=\bar{A}$;
\item $A$ preserves the form $\omega$;
\item $\|A\|_2<\infty$.
\item $A$ is invariant on $H$, i.e., $A(H)\subseteq H$.
\end{enumerate}
\end{remark}

\begin{remark}\label{def.SPInf2}
By Corollary \ref{AASharp}, the definition of $\Sp(\infty)$ is also
equivalent to
\begin{enumerate}
\item $A=\bar{A}$;
\item $A^\#A=AA^\#=id$;
\item $\|A\|_2<\infty$.
\end{enumerate}
\end{remark}

\begin{proposition}
$\Sp(\infty)$ is a group.
\end{proposition}

\begin{proof}
First we show that if $A\in \Sp(\infty)$, then $A^{-1} \in
\Sp(\infty)$. By the assumption on $A$, it is easy to verify that
$A^{-1}$ satisfies (1), (2), (3) and (5) in Remark \ref{def.SPInf1}.
We need to show that $A^{-1}$ satisfies the condition (4), i.e.
$\|A^{-1}\|_2<\infty$. Suppose
\[
A=\left(
\begin{array}{ll}
a &b\\
\bar{b} &\bar{a}
\end{array}
\right) \hspace{.1in}\mbox{and}\hspace{.1in} A^{-1}=\left(
\begin{array}{ll}
a' &b'\\
\overline{b'} &\overline{a'}
\end{array}
\right),
\]
where by our assumptions all blocks are bounded operators, and in
addition $b$ is a Hilbert-Schmidt operator. We want to prove $b'$ is
also a Hilbert-Schmidt operator. $AA^{-1}=I$ and $A^{-1}A=I$ imply
that
\[
ab^{\prime}=-b\overline{a'}, \hskip0.1in  a'a+b'\bar{b}=I.
\]
The last equation gives $a'ab'+b'\bar{b}b'=b'$, and so

\[
b'=a'ab'+b'\bar{b}b'=-a'b\overline{a'}+b'\bar{b}b'
\]
which is a Hilbert-Schmidt operator as $b$ and $\bar{b}$ are
Hilbert-Schmidt. Therefore $\|A^{-1}\|_2<\infty$ and $A^{-1}\in
\Sp(\infty)$.

Next we show that if $A,B\in \Sp(\infty)$, then $AB\in \Sp(\infty)$.
By the assumption on $A$ and $B$, it is easy to verify that $AB$
satisfies (1), (2), (3) and (5) in Remark \ref{def.SPInf1}. We need
to show that $AB$ satisfies the condition (4), i.e.
$\|AB\|_2<\infty$. Suppose
\[
A=\left(
\begin{array}{ll}
a &b\\
\bar{b} &\bar{a}
\end{array}
\right) \hspace{.1in}\mbox{and}\hspace{.1in} B=\left(
\begin{array}{ll}
c &d\\
\bar{d} &\bar{c}
\end{array}
\right),
\]
where all blocks are bounded, and $\Vert b\Vert_{HS}, \Vert
d\Vert_{HS}<\infty$. Then
\[
AB=\left(
\begin{array}{ll}
ac+b\bar{d} & ad+b\bar{c}\\
\bar{b}c+\bar{a}\bar{d} & \bar{b}d+\bar{a}\bar{c}
\end{array}
\right).
\]
Then
\[
\Vert AB \Vert_{2}^{2}=\|ad+b\bar{c}\|_{HS} \leqslant
\|ad\|_{2}+\|b\bar{c}\|_{HS}<\infty,
\] since both $ad$ and
$b\bar{c}$ are Hilbert-Schmidt operators. Therefore
$\|AB\|_2<\infty$ and $AB\in \Sp(\infty)$.
\end{proof}

\section{Symplectic Representation of $\Diff(S^1)$}

\begin{definition}\label{d.4.1}
Let $\Diff(S^1)$ be the group of orientation preserving $C^\infty$
diffeomorphisms of $S^1$. $\Diff(S^1)$ acts on $H$ as follows
\[
(\phi.u)(\theta)=u(\phi^{-1}(\theta))
-\frac{1}{2\pi}\int_0^{2\pi} u(\phi^{-1}(\theta)) d\theta.
\]
\end{definition}
Note that if $u\in H$ is real-valued, then $\phi.u$ is real-valued
as well.

\begin{proposition}
The action of $\Diff(S^1)$ on $H$ gives a group homomorphism
\[
\Phi:\Diff(S^1)\to \Aut H
\]
defined by $\Phi(\phi)(u)=\phi.u$, for $\phi\in\Diff(S^1)$ and $u\in H$,
where $\Aut H$ is the group of automorphisms on $H$.
\end{proposition}

\begin{proof}
Let $u\in H$, then $\phi.u$ is  a $C^{\infty}$ function with the
mean value $0$, and so $\phi.u\in H$. It is also clear that
$\phi.(u+v)=\phi.u+\phi.v$ and $\phi.(\lambda u)=\lambda\phi.u$. So
$\Phi$ is well--defined as a map from $\Diff(S^1)$ to $\End H$, the
space of endomorphisms on $H$. Now let us check that $\Phi$ is a
group homomorphism. Suppose $\phi,\psi\in\Diff(S^1)$ and $u\in H$,
then
\begin{align*}
&
\Phi(\phi\psi)(u)(\theta)
=u\big((\phi\psi)^{-1}(\theta)\big)
-\frac{1}{2\pi}\int_0^{2\pi} u\big((\phi\psi)^{-1}(\theta)\big)d\theta
\\&
=u\big((\psi^{-1}\phi^{-1})(\theta)\big)
-\frac{1}{2\pi}\int_0^{2\pi} u\big((\psi^{-1}\phi^{-1})(\theta)\big)d\theta.
\end{align*}
On the other hand,
\begin{align*}
& \Phi(\phi)\Phi(\psi)(u)(\theta)
=\Phi(\phi)\left[u(\psi^{-1}(\theta)) -\frac{1}{2\pi}\int_0^{2\pi}
u(\psi^{-1}(\theta))d\theta\right]
\\&
=\Phi(\phi)\left[u(\psi^{-1}(\theta))\right]
=u\big((\psi^{-1}\phi^{-1})(\theta)\big)
-\frac{1}{2\pi}\int_0^{2\pi}
u\big((\psi^{-1}\phi^{-1})(\theta)\big)d\theta.
\end{align*}
So $\Phi(\phi\psi)=\Phi(\phi)\Phi(\psi)$. In particular, the image
of $\Phi$ is in the $\Aut H$.
\end{proof}

\begin{lemma}\label{l.4.4}
Any $\phi\in\Diff(S^1)$  preserves the form $\omega$, that is,
$\omega(\phi.u,\phi.v)=\omega(u,v)$ for any $ u,v\in H$.
\end{lemma}

\begin{proof}
By Definition \ref{d.4.1} $\phi.u=u(\psi)-u_0,\phi.v=v(\psi)-v_0$,
where $\psi=\phi^{-1}$ and $u_0,v_0$ are the constants. Then
\begin{align*}
\omega(\phi.u,\phi.v)
&=\omega(u(\psi)-u_0,v(\psi)-v_0)
\\&
=\frac{1}{2\pi}\int_0^{2\pi}
\big(u(\psi(\theta))-u_0\big)
\big(v(\psi(\theta))-v_0\big)' d\theta
\\&
=\frac{1}{2\pi}\int_0^{2\pi}
u(\psi)v'(\psi)\psi'(\theta)d\theta
-\frac{1}{2\pi}\int_0^{2\pi}
u_0v(\psi(\theta))d\theta
\\&
=\frac{1}{2\pi}\int_0^{2\pi}u(\psi)v'(\psi)d\psi
\\&
=\omega(u,v).
\end{align*}
\end{proof}

We are going to prove that a diffeomorphism $\phi\in\Diff(S^1)$ acts
on $H$ as a bounded linear map, and that $\Phi(\phi)$ is in
$\Sp(\infty)$. The next lemma is a generalization of a proposition
in a paper of G.~Segal\cite{Segal1981}.

\begin{lemma}\label{lem.BigLemma}
Let $\psi\neq id \in\Diff(S^1)$ and $\phi=\psi^{-1}$.
Let
\[
I_{n,m}=(\psi.e^{im\theta},e^{in\theta})
=\frac{1}{2\pi}\int_0^{2\pi}e^{im\phi-in\theta}d\theta.
\]
Then
\begin{enumerate}
\item
$\displaystyle \sum_{n>0,m<0}|n||I_{n,m}|^2 <\infty$, and
$\displaystyle \sum_{m>0,n<0}|n||I_{n,m}|^2 <\infty$.

\item
For sufficiently large $|m|$ there is a constant $C$ independent of
$m$ such that
\begin{equation}\label{e.4.1}
\sum_{n\neq0}|n||I_{n,m}|^2 < C|m|.
\end{equation}
\end{enumerate}
\end{lemma}

\begin{proof}
Let
\[
m_{\phi'}=\min\{\phi'(\theta)|\theta\in S^1\};
\hspace{.1in}\mbox{and}\hspace{.1in}
M_{\phi'}=\max\{\phi'(\theta)|\theta\in S^1\}.
\]
Since $\phi$ is a diffeomorphism, we have $0<m_{\phi'}<M_{\phi'}<\infty$.

Take four points $a,b,c,d$ on the unit circle such that
$a$ corresponds to $m_{\phi'}$ in the sense $\tan(a)=m_{\phi'}$,
$b$ corresponds to $M_{\phi'}$ in the sense $\tan(b)=M_{\phi'}$,
$c$ is opposite to $a$, i.e., $c=a+\pi$,
$d$ is opposite to $b$, i.e., $d=b+\pi$.
The four points on the circle are arranged in the counter-clockwise order,
and $0<a<b<\frac{\pi}{2}$, $\pi<c<d<\frac{3}{2}\pi$.

Let $\tau\in S^1$ such that $\tau\neq\frac{\pi}{4},\frac{5}{4}\pi$.
Define a function $\phi_\tau$ on $S^1$ by
\[
\phi_\tau(\theta)=
\frac{\cos\tau\cdot\phi(\theta)-\sin\tau\cdot\theta}
{\cos\tau-\sin\tau}.
\]

We will show that if $\tau\in (b,c)$ or $\tau\in(d,a)$,
then $\phi_\tau$ is an orientation preserving diffeomorphism of $S^1$,
where $(b,c)$ is the open arc from the point $b$ to the point $c$,
and $(d,a)$ is the open arc from the point $d$ to the point $a$.

Clearly $\phi_\tau$ is a $C^\infty$ function on $S^1$.
Also, $\phi_\tau(0)=0$ and $\phi_\tau(2\pi)=2\pi$.
Taking derivative with respect to $\theta$, we have
\[
\phi'_\tau(\theta)=
\frac{\cos\tau\cdot\phi'(\theta)-\sin\tau}
{\cos\tau-\sin\tau}.
\]
By the choice of $\tau$, we can prove that $\phi'_\tau(\theta)>0$.
Therefore, $\phi_\tau$ is an orientation preserving diffeomorphism as claimed.

Let $m,n\in\mathbb{Z}\backslash\{0\}$.
Let $\tau_{mn}=\Arg(m+in)$, i.e., the argument of the complex number
$m+in$, considered to be in $[0,2\pi]$.
Then we have $m\phi-n\theta=(m-n)\phi_{\tau_{mn}}$.

If $\tau_{mn}\in(b,c)$, then $\phi_{\tau_{mn}}$ is a diffeomorphism.
Let $\psi_{\tau_{mn}}=\phi_{\tau_{mn}}^{-1}$. Then
\[
I_{n,m}
=\frac{1}{2\pi}\int_0^{2\pi}e^{i(m-n)\phi_{\tau_{mn}}}d\theta
=\frac{1}{2\pi}\int_0^{2\pi}e^{i(m-n)\theta}
\psi_{\tau_{mn}}'(\theta)d\theta,
\]
where the last equality is by change of variable.
On integration by parts $k$ times, we have
\[
I_{n,m}
=\left(\frac{1}{i(m-n)}\right)^k\frac{1}{2\pi}\int_0^{2\pi}
e^{i(m-n)\theta}\psi_{\tau_{mn}}^{(k+1)}(\theta)d\theta.
\]

Let $\alpha=[\alpha_0,\alpha_1]$ be a closed arc contained in the arc $(b,c)$.
Let $S_\alpha$ be the set of all pairs of nonzero integers $(m,n)$ such that
$\alpha_0<\tau_{mn}<\alpha_1$, where $\tau_{mn}=\Arg(m+in)$.
We are going to consider an upper bound of the sum
$\sum_{(m,n)\in S_\alpha}|n||I_{n,m}|^2$.

For the pair $(m,n)$, if $|m-n|=p$, the condition
$\alpha_0<\tau_{mn}<\alpha_1$ gives us both an upper bound and a lower
bound for $n$:
\[
\frac{m_{\phi'}}{m_{\phi'}-1}p
\le n \le
\frac{M_{\phi'}}{M_{\phi'}-1}p.
\]
So $|n|\le C_1p$ where $C_1$ is a constant which does not depend on
the pair $(m,n)$. Also, the number of pairs $(m,n)\in S_\alpha$ such
that $|m-n|=p$ is bounded by $C_2p$ for some constant $C_2$. Let
$C_3=\max \Big\{ |\psi_\tau^{(k+1)}(\theta)| : \theta\in S^1,
\tau\in[\alpha_0,\alpha_1] \Big\}$. Then
\[
|I_{n,m}|\le
C_3\Big|\frac{1}{i(m-n)}\Big|^k\frac{1}{2\pi}\int_0^{2\pi}
e^{i(m-n)\theta}d\theta =C_3 p^{-k}.
\]
Therefore,
\begin{align*}
\sum_{(m,n)\in S}|n||I_{n,m}|^2
&=\sum_p\sum_{(m,n)\in S_\alpha; |m-n|=p}|n||I_{n,m}|^2
\\&
\le\sum_p C_1p\cdot C_3^2p^{-2k}\cdot C_2p =C_\alpha  \sum_p
p^{-(2k-2)},
\end{align*}
where the constant $C_\alpha$ depends on the arc $\alpha$.

Similarly, for a closed arc $\beta=[\beta_0,\beta_1]$ contained in
the arc $(d,a)$, we have
\[
\sum_{(m,n)\in S_\beta}|n||I_{n,m}|^2 \le C_\beta \sum_p
p^{-(2k-2)},
\]
where the constant $C_\beta$ depends on the arc $\beta$.

Now let $\alpha=[\frac{\pi}{2},\pi]$, and $\beta=[\frac{3}{2}\pi,2\pi]$.
Then $\alpha$ is contained in $(b,c)$ and $\beta$ is contained in $(d,a)$.
We have
\[
\sum_{n>0,m<0}|n||I_{n,m}|^2
=C_\alpha \cdot \sum_p p^{-(2k-2)}
<\infty
\]
and
\[
\sum_{n<0,m>0}|n||I_{n,m}|^2 =C_\beta \cdot \sum_p p^{-(2k-2)}
<\infty,
\]
which proves (1) of the lemma.

To prove (2), we let $\alpha=[\alpha_0,\alpha_1]$ be a closed arc
contained in the arc $(b,c)$ such that $b<\alpha_0<\frac{\pi}{2}$
and $\pi<\alpha_1<c$, and $\beta=[\beta_0,\beta_1]$ be a closed arc
contained in the arc $(d,a)$ such that $d<\beta_0<\frac{3}{2}\pi$
and $0<\beta_1<a$. Then we have
\[
\sum_{(m,n)\in S_\alpha}|n||I_{n,m}|^2 +\sum_{(m,n)\in
S_\beta}|n||I_{n,m}|^2 \leqslant  C_{\alpha\beta}
\]
for some constant $C_{\alpha\beta}$.

Let $m>0$ be sufficiently large, and $N_m$ be the largest integer
less than or equal to $m\tan(\alpha_0)$,
\[
\sum_{0<n\leqslant N_m}|I_{n,m}|^2 \leqslant
\sum_{n\neq0}|I_{n,m}|^2.
\]

Note that $I_{n,m}$ is the $n$th Fourier coefficient of
$\psi.e^{im\theta}$. Therefore,
\[
\sum_{n\neq0}|I_{n,m}|^2
=\|\psi.e^{im\theta}\|_{L^2}
\]
which is bounded by a constant $K$.
Therefore,
\[
\sum_{0<n\leqslant N_m}|n||I_{n,m}|^2 \leqslant
Km\tan\left(\alpha_0\right).
\]
On the other hand,
\[
\sum_{n<0}|n||I_{n,m}|^2 + \sum_{n>N_m}|n||I_{n,m}|^2 \leqslant
\sum_{(m,n)\in S_\alpha}|n||I_{n,m}|^2 + \sum_{(m,n)\in
S_\beta}|n||I_{n,m}|^2 =C_{\alpha\beta}.
\]
Therefore,
\[
\sum_{n\neq 0}|n||I_{n,m}|^2 \leqslant C_{\alpha\beta} +
Km\tan(\alpha_0) \leqslant m C_+,
\]
where $C_+$ can be chosen to be, for example,
$K\tan(\alpha_0)+C_{\alpha\beta}$, which is independent of $m$.

Similarly, for $m<0$ with sufficiently large $|m|$
\[
\sum_{n\neq 0}|n||I_{n,m}|^2 \leqslant m C_-.
\]

Let $C=\max \{C_+,C_-\}$. Then we have, for sufficiently
large $|m|$,
\[
\sum_{n\neq 0}|n||I_{n,m}|^2 \leqslant |m|C,
\]
which proves (2) of the lemma.
\end{proof}

\begin{lemma}\label{l.4.6}
For any $\psi\in\Diff(S^1)$, $\Phi(\psi)\in B(H)$, the space of
bounded linear maps on $H$. Moreover,
\[
\Vert\Phi(\psi)\Vert \leqslant C, \
\Vert\Phi(\psi)\Vert_{2}\leqslant C, \] where $C$ is the constant in
Equation \ref{e.4.1}.
\end{lemma}

\begin{proof}
First observe that the operator norm of $\Phi(\psi)$ is
\[
\|\Phi(\psi)\|=\sup\{\|\psi.u\|_\omega \hspace{.2cm}|\hspace{.2cm}
u\in H,\|u\|_\omega=1\}.
\]
For any $u\in H$, let $\hat{u}$ be its Fourier coefficients, that is
$\hat{u}(n)=(u,\hat{e}_n)$, and let $\tilde{u}$ be defined by
$\tilde{u}=(u,\tilde{e}_n)_\omega$
(\ref{nota.Inner},\ref{def.OmegaBasis}). It can be verified that the
relation between $\hat{u}$ and $\tilde{u}$ is: if $n>0$, then
$\tilde{u}(n)=\sqrt{n}\hat{u}(n)$; if $n<0$, then
$\tilde{u}(n)=i\sqrt{|n|}\hat{u}(n)$. We have
\[
\|u\|_\omega^2
=(u,u)_\omega
=(\tilde{u},\tilde{u})_{l^2}
=\sum_{n\neq0}|\tilde{u}(n)|^2
=\sum_{n\neq0}|n||\hat{u}(n)|^2.
\]
Let $\phi=\psi^{-1}$. We have $u(\phi)=\sum_{m\neq0}\hat{u}(m)e^{im\phi}$.
Using the notation $I_{n,m}$ (\ref{lem.BigLemma}), we have
\begin{align*}
& \|\psi.u\|_\omega^2 =\sum_{n\neq0}|n||\widehat{\psi.u}(n)|^2
=\sum_{n\neq0}|n|\Big| \frac{1}{2\pi}\int_0^{2\pi}
u(\phi(\theta))e^{-in\theta}d\theta \Big|^2
\\
&
=\sum_{n\neq0}|n|\Big| \frac{1}{2\pi}\int_0^{2\pi}
\sum_{m\neq0}\hat{u}(m)e^{im\phi} e^{-in\theta}d\theta \Big|^2
\\&
=\sum_{n\neq0}|n|\Big|\sum_{m\neq0}\hat{u}(m)
\frac{1}{2\pi}\int_0^{2\pi}e^{im\phi-in\theta}d\theta \Big|^2
\\&
=\sum_{n\neq0}|n|\Big|\sum_{m\neq0}\hat{u}(m) I_{n,m}\Big|^2 \\
& \leqslant \sum_{m,n\neq0}|n||\hat{u}(m)|^2|I_{n,m}|^2
=\sum_{m\neq0}|\hat{u}(m)|^2 \sum_{n\neq0}|n||I_{n,m}|^2
\\&
=\sum_{|m|\leqslant M_0}|\hat{u}(m)|^2 \sum_{n\neq0}|n||I_{n,m}|^2
+\sum_{|m|>M_0}|\hat{u}(m)|^2 \sum_{n\neq0}|n||I_{n,m}|^2,
\end{align*}
where the constant $M_0$ in the last equality is a positive integer
large enough so that we can apply part (2) of Lemma
\ref{lem.BigLemma}. It is easy to see that the first term in the
last equality is finite. For the second term we use Lemma
\ref{lem.BigLemma}
\[
\sum_{|m|>M_0}|\hat{u}(m)|^2 \sum_{n\neq0}|n||I_{n,m}|^2 \leqslant
C\sum_{|m|>M_0}|\hat{u}(m)|^2 |m| \leqslant C.
\]
Thus  for any  $u\in H$ with $\|u\|_\omega=1$, $\|\psi.u\|_\omega$
is uniformly bounded. Therefore, $\Phi(\psi)$ is a bounded operator
on $H$.

Now we can use Lemma \ref{lem.BigLemma} again to estimate the norm
$\Vert \Phi(\psi) \Vert_{2}$

\begin{align*}
& \|\Phi(\psi)\|_{2}
=\sum_{n>0,m<0}|(\psi.\tilde{e}_m,\tilde{e}_n)_\omega|^2
=\sum_{n>0,m<0}|n||(\psi.\hat{e}_m,\hat{e}_n)|^2
\\&
=\sum_{n>0,m<0}|n||I_{n,m}|^2 <\infty.
\end{align*}
\end{proof}

\begin{theorem}\label{t.4.6}
$\Phi:\Diff(S^1)\to \Sp(\infty)$ is a group homomorphism. Moreover,
$\Phi$ is injective, but not surjective.
\end{theorem}

\begin{proof}
Combining Lemma \ref{l.4.4} and Lemma \ref{l.4.6} we see that for
any diffeomorphism $\psi\in\Diff(S^1)$ the map $\Phi(\psi)$ is an
invertible bounded operator on $H$, it preserves the form $\omega$,
and $\|\Phi(\psi)\|_{2}<\infty$. In addition, by our remark after
Definition \ref{d.4.1} $\psi.u$ is real-valued, if $u$ is
real-valued. Therefore, $\Phi$ maps $\Diff(S^1)$ into $\Sp(\infty)$.

Next, we first prove that $\Phi$ is injective. Let
$\psi_1,\psi_2\in\Diff(S^1)$, and denote $\phi_1=\psi_1^{-1},
\phi_2=\psi_2^{-1}$. Suppose $\Phi(\psi_1)=\Phi(\psi_2)$, i.e.
$\psi_1.u=\psi_2.u$, for any $u\in H$. In particular,
$\psi_1.e^{i\theta}=\psi_2.e^{i\theta}$. Therefore
\[
e^{i\phi_1}-C_1=e^{i\phi_2}-C_2,
\]
where $C_1=\frac{1}{2\pi}\int_0^{2\pi}e^{i\phi_1}d\theta$, and
$C_2=\frac{1}{2\pi}\int_0^{2\pi}e^{i\phi_2}d\theta$. Note that
$e^{i\phi_1}$ and $e^{i\phi_2}$ have the same image as maps from
$S^{1}$ to $\mathbb{C}$. This implies $C_1=C_2$, since otherwise
$e^{i\phi_1}=e^{i\phi_2}+(C_1-C_2)$ and $e^{i\phi_1}$ and
$e^{i\phi_2}$ would have had different images. Therefore, we have
$e^{i\phi_1}=e^{i\phi_2}$. But the function $e^{i\tau}:S^1\to S^1$
is an injective function, so $\phi_1=\phi_2$. Therefore
$\psi_1=\psi_2$, and so $\Phi$ is injective.

To prove that $\Phi$ is not surjective, we will construct an
operator $A\in \Sp(\infty)$ which can not be written as $\Phi(\psi)$
for any $\psi\in\Diff(S^1)$. Let the linear map $A$ be defined by
the corresponding matrix $\{ A_{m, n}\}_{m, n \in \mathbb{Z}}$ with
the entries
\begin{align*}
& A_{1,1}=A_{-1,-1}=\sqrt{2}\\
& A_{1,-1}=i, A_{-1,1}=-i\\
& A_{m,m}=1, \hspace{.1in}\mbox{for }m \neq \pm1
\end{align*}
with all other entries  being $0$.

First we show that $A\in \Sp(\infty)$. For any $u\in H$, we can
write $u=\sum_{n\neq0}\tilde{u}(n)\tilde{e}_n$. Then $A$ acting on
$u$ changes  only $\tilde{e}_1$ and $\tilde{e}_{-1}$ . Therefore,
$Au\in H$, and clearly $A$ is a well--defined  bounded linear map on
$H$ to $H$. Moreover, $\|A\|_{2}<\infty$.  It is clear that
$A_{m,n}=\overline{A_{-m,-n}}$, and therefore $A=\bar{A}$ by
Proposition \ref{SomeFacts1}. Moreover, $A$ preserves the form
$\omega$ by  part(II) of Proposition \ref{prop.Preserve}, as

\[\sum_{k\neq 0}\sgn(mk)A_{k,m}\overline{A_{k,n}}
=\delta_{m,n}.
\]
Finally, $A$ is invertible, since $\{A_{k, m}\}_{m, n \in
\mathbb{Z}}$ is, with the inverse $\{B_{k, m}\}_{m, n \in
\mathbb{Z}}$ given by
\begin{align*}
& B_{1,1}=B_{-1,-1}=\sqrt{2}\\
& B_{1,-1}=-i, B_{-1,1}=i\\
& B_{m,m}=1, \hspace{.1in}\mbox{for }m \neq \pm1
\end{align*}
with all other entries being $0$. Next we show that
$A\not=\Phi(\psi)$ for any $\psi\in\Diff(S^1)$. First observe that
if we look at any basis element $\tilde{e}_1=e^{i\theta}$ as a
function from $S^1$ to $\mathbb{C}$, then the image of this function
lies on the unit circle. Clearly, when acted by a diffeomorphism
$\phi\in\Diff(S^1)$, the image of the function $\phi.e^{i\theta}$ is
still a circle with radius $1$. But if we consider $A\tilde{e}_1$ as
a function from $S^1$ to $\mathbb{C}$, we will show that the image
of the function $A\tilde{e}_1:S^1\to\mathbb{C}$ is not  a circle.
Therefore, $A\not=\Phi(\psi)$ for any $\psi\in\Diff(S^1)$. Indeed,
by definition of $A$ we have
\[
A\tilde{e}_1=\sqrt{2}\tilde{e}_1 - i \tilde{e}_{-1}.
\]
Let us write it as a function on $S^1$
\[
A\tilde{e}_1(\theta)=\sqrt{2}e^{i\theta}-e^{-i\theta}
=(\sqrt{2}-1)\cos\theta+i(\sqrt{2}+1)\sin\theta,
\]
and then we see that the image lies on an ellipse, which is not the
unit circle
\[ \frac{x^2}{(\sqrt{2}-1)^2}+\frac{y^2}{(\sqrt{2}+1)^2}=1. \]
\end{proof}

\section{The Lie algebra associated with $\Diff(S^1)$}

Let
$\diff(S^1)$ be the space of smooth vector fields on $S^1$. Elements
in $\diff(S^1)$ can be identified with smooth functions on $S^1$.
The space $\diff(S^1)$ is a Lie algebra with the following Lie
bracket
\[
[X,Y]=XY'-X'Y, \hskip0.1in X,Y\in\diff(S^1),
\]
where $X'$ and $Y'$ are derivatives with respect to $\theta$.

Let $X\in\diff(S^1)$, and $\rho_t$ be the corresponding flow of
diffeomorphisms. We define an action of $\diff(S^1)$ on $H$ as
follows: for $X\in\diff(S^1)$ and $u\in H$, $X.u$ is a function on
$S^1$ defined by
\[
(X.u)(\theta)=\left.\frac{d}{dt}\right|_{t=0} \left[
(\rho_t.u)(\theta) \right],
\]
where $\rho_t$ acts on $u$ via the representation
$\Phi:\Diff(S^1)\to \Sp(\infty)$.

The next proposition shows that the action is well--defined, and
also gives an explicit formula of $X.u$.

\begin{proposition}\label{prop.diffAction}
Let $X\in\diff(S^1)$.
Then
\[
(X.u)(\theta)
=u'(\theta)(-X(\theta))
-\frac{1}{2\pi}\int_0^{2\pi}u'(\theta)(-X(\theta))d\theta,
\]
that is, $X.u$ is the function $-u'X$ with the $0$th Fourier
coefficient replaced by $0$.
\end{proposition}

\begin{proof}
Let $\rho_t$ be the flow that corresponds to $X$, and $\lambda_t$ be the flow
that corresponds to $-X$. Then $\lambda_t$ is the inverse of $\rho_t$
for all $t$.
\[
(X.u)(\theta) =\left.\frac{d}{dt}\right|_{t=0} \left[
(\rho_t.u)(\theta) \right]
=\left.\frac{d}{dt}\right|_{t=0}\bigg[u(\lambda_t(\theta))
-\frac{1}{2\pi}\int_0^{2\pi}u(\lambda_t(\theta))d\theta\bigg].
\]
Using the chain rule, we have
\[
\left.\frac{d}{dt}\right|_{t=0}u(\lambda_t(\theta))
=u'(\theta)(-\widetilde{X}(\theta)),
\]
and
\[
\left.\frac{d}{dt}\right|_{t=0}\frac{1}{2\pi}
\int_0^{2\pi}u(\lambda_t(\theta))d\theta
=\frac{1}{2\pi}\int_0^{2\pi}u'(\theta)(-X(\theta))d\theta.
\]
\end{proof}

\begin{notation}
We consider $\diff(S^1)$ as a subspace of the space of real-valued $L^2$ functions on $S^1$.
The space of real-valued $L^2$ functions on $S^1$ has an orthonormal basis
\[
\mathcal{B}=\{X_{l}=\cos(m\theta),  Y_{k}=\sin(k\theta), l=0, 1,
..., k=1, 2,...\}
\]
which is contained in $\diff(S^1)$.
\end{notation}
Let us consider how these basis elements act on $H$.

\begin{proposition}
For  any $l=0, 1, ..., k=1, 2,...$ the basis elements $X_l,Y_k$ act
on $H$ as linear maps. In the basis $\mathcal{B}_\omega$ of $H$,
they are represented by infinite dimensional matrices with $(m,n)$th
entries equal to

\begin{align*}
(X_l)_{m,n}&=(X_l.\tilde{e}_n,\tilde{e}_m)_\omega=
s(m,n)\frac{1}{2}\sqrt{|mn|}(\delta_{m-n,l}+\delta_{n-m,l})\\
(Y_k)_{m,n}&=(Y_k.\tilde{e}_n,\tilde{e}_m)_\omega=
s(m,n)(-i)\frac{1}{2}\sqrt{|mn|}(\delta_{m-n,k}-\delta_{n-m,k})
\end{align*} where $m,n\neq0$,
\[
s(m,n)=\left\{
\begin{array}{ll}
-i & m,n>0\\
1  & m>0,n<0\\
1  & m<0,n>0\\
i  & m,n<0.
\end{array}
\right.
\]

\end{proposition}

\begin{proof}
By Proposition \ref{prop.diffAction} and a simple verification
depending on the signs of $m,n$ we see that
\begin{align*}
X_l.e^{in\theta}&=-ine^{in\theta}\cos(l\theta)
=-\frac{1}{2}in\left[e^{i(n+l)\theta}+e^{i(n-l)\theta}\right]\\
Y_k.e^{in\theta}&=-ine^{in\theta}\sin(k\theta)
=-\frac{1}{2}n\left[e^{i(n+k)\theta}-e^{i(n-k)\theta}\right].
\end{align*}
Indeed, recall that a basis element
$\tilde{e}_n\in\mathcal{B}_\omega$ has the form
\[
\tilde{e}_n=\left\{
\begin{array}{ll}
\frac{1}{\sqrt{n}}e^{in\theta} & n>0\\
\frac{1}{i\sqrt{|n|}}e^{in\theta} & n<0.
\end{array}
\right.
\]
Suppose $m,n>0$
\[
X_l.\tilde{e}_n=\frac{1}{\sqrt{n}}X_l.e^{in\theta}
=-\frac{1}{2}i\sqrt{n}\left[e^{i(n+l)\theta}+e^{i(n-l)\theta}\right],
\]
and
\[
(e^{i(n+l)\theta},\tilde{e}_m)_\omega=\sqrt{m}\delta_{m-n,k};
\hspace{.2in}
(e^{i(n-l)\theta},\tilde{e}_m)_\omega=\sqrt{m}\delta_{n-m,l}.
\]
Therefore,
\[
(X_l)_{m,n}=(X_l.\tilde{e}_n,\tilde{e}_m)_\omega=
(-i)\frac{1}{2}\sqrt{|mn|}(\delta_{m-n,l}+\delta_{n-m,l}).\\
\]
All other cases can be verified similarly.
\end{proof}

\begin{remark}
Recall that $\mathbb{H}_\omega$ is the completion of $H$ under
the metric $(\cdot,\cdot)_\omega$.
The above calculation shows that the trigonometric basis
$X_l,Y_k$ of $\diff(S^1)$ act on $\mathbb{H}_\omega$
as \emph{unbounded} operators. They are densely defined on the subspace
$H\subseteq\mathbb{H}_\omega$.
\end{remark}

\section{Brownian motion on $\Sp(\infty)$}\label{s.6}

\begin{notation}
As in \cite{AirMall2001}, let $\mathfrak{sp}\left(\infty\right)$  be
the set of infinite-dimensional matrices $A$  which can be written
as block matrices of the form
\[
\left(
\begin{array}{ll}
a &b\\
\bar{b} &\bar{a}
\end{array}
\right)
\]
such that $a+a^\dag=0$, $b=b^T$, and $b$ is a Hilbert-Schmidt
operator.
\end{notation}

\begin{remark}
 The set $\mathfrak{sp}\left(\infty\right)$ has a structure
of Lie algebra with the operator commutator as a Lie bracket, and we
associate this Lie algebra with the group $\Sp(\infty)$.
\end{remark}




\begin{proposition}\label{GEntry}
Let $\{A_{m,n}\}_{m, n \in \mathbb{Z}\backslash\{0\}}$ be the matrix
corresponding to an operator $A$. Then any $A\in
\mathfrak{sp}\left(\infty\right)$ satisfies (1)
$A_{m,n}=\overline{A_{-m,-n}}$; (2) $A_{m,n}+\overline{A_{n,m}}=0$,
for $m,n>0$; (3) $A_{m,n}=A_{-n,-m}$, for $m>0,n<0$.

Moreover, $A\in \mathfrak{sp}\left(\infty\right)$ if and only if (1)
$A=\bar{A}$; (2) $\pi^+A\pi^-$ is Hilbert-Schmidt; (3) $A+A^\#=0$.

\end{proposition}

\begin{proof}
The first part follows directly from definition of
$\mathfrak{sp}\left(\infty\right)$.  Then we can use this fact and
the formula for the matrix entries of $A^\#$ in Proposition
\ref{SomeFacts1} to prove the second part.
\end{proof}


\begin{definition}
Let $HS$ be the space of Hilbert-Schmidt matrices viewed as a real
vector space, and $\SpHS=\mathfrak{sp}\left(\infty\right) \cap HS$.
\end{definition}

The space $HS$ as a real Hilbert space has an orthonormal basis
\[
\mathcal{B}_{HS}=\{e_{mn}^{Re}:m,n\neq0\}\cup\{e_{mn}^{Im}:m,n\neq0\},
\]
where $e_{mn}^{Re}$ is a matrix with $(m,n)$-th entry $1$ all other
entries $0$, and $e_{mn}^{Im}$ is a matrix with $(m,n)$th entry $i$
all other entries $0$.

The space $\SpHS$ is a closed subspace of $HS$, and therefore
 a real Hilbert space. According to the symmetry of the
matrices in $\SpHS$, we define a projection $\pi:HS\to \SpHS$, such
that
\begin{align*}\label{Basis}
\pi(e_{mn}^{Re})
&=\frac{1}{2}\big(e_{mn}^{Re}-e_{nm}^{Re}+e_{-m,-n}^{Re}-e_{-n,-m}^{Re}\big),
&\mbox{if } \sgn(mn)>0\\
\pi(e_{mn}^{Im})
&=\frac{1}{2}\big(e_{mn}^{Im}+e_{nm}^{Im}-e_{-m,-n}^{Im}-e_{-n,-m}^{Im}\big),
&\mbox{if } \sgn(mn)>0\\
\pi(e_{mn}^{Re})
&=\frac{1}{2}\big(e_{mn}^{Re}+e_{-n,-m}^{Re}+e_{-m,-n}^{Re}+e_{n,m}^{Re}\big),
&\mbox{if } \sgn(mn)<0\\
\pi(e_{mn}^{Im})
&=\frac{1}{2}\big(e_{mn}^{Im}+e_{-n,-m}^{Im}-e_{-m,-n}^{Im}-e_{nm}^{Im}\big),
&\mbox{if } \sgn(mn)<0
\end{align*}

\begin{notation}
We choose $\mathcal{B}_{\SpHS}=\pi(\mathcal{B}_{HS})$ to be the
orthonormal basis of $\SpHS$.
\end{notation}

Clearly, if $A\in \SpHS$, then $|A|_{\SpHS}=|A|_{HS}$.

\begin{definition}
Let $W_{t}$ be a Brownian motion on $\SpHS$ which has the mean zero
and covariance $Q$, where $Q$ is assumed to be a positive symmetric
trace class operator on $H$. We further assume that $Q$ is diagonal
in the basis $\mathcal{B}_{\SpHS}$.
\end{definition}

\begin{remark}
$Q$ can also be viewed as a positive function on the set
$\mathcal{B}_{\SpHS}$, and the Brownian motion $W_{t}$ can be
written as
\begin{equation}
W_{t}=\sum_{\xi\in\mathcal{B}_{\SpHS}}\sqrt{Q(\xi)}B_{t}^\xi\xi,
\end{equation}
where $\{B_{t}^\xi\}_{\xi\in\mathcal{B}_{\SpHS}}$ are standard
real-valued mutually independent Brownian motions.
\end{remark}

Our goal now is to construct a Brownian motion on the group
$\Sp(\infty)$ using the Brownian motion $W_{t}$ on $\SpHS$. This is
done by solving the Stratonovich stochastic differential equation
\begin{equation}
\delta X_{t} = X_{t} \delta W_{t}.
\end{equation}
This equation can be written as the following It\^o stochastic
differential equation
\begin{equation}\label{ItoSDE}
dX_{t} = X_{t}dW_{t} + \frac{1}{2}X_{t}Ddt,
\end{equation}
where $D=\Diag(D_m)$ is a diagonal matrix with entries
\begin{equation}\label{DMatrix}
D_m=-\frac{1}{4}\sgn(m)\sum_k\sgn(k)\left[Q_{mk}^{Re}+Q_{mk}^{Im}\right]
\end{equation}
with $Q_{mk}^{Re}=Q(\pi(e_{mk}^{Re}))$ and $Q_{mk}^{Im}=Q(\pi(e_{mk}^{Im}))$.

\begin{notation}
Denote by $\SpQHS=Q^{1/2}(\SpHS)$ which is a subspace of $\SpHS$.
Define an inner product on $\SpQHS$ by $\langle u,v \rangle_{\SpQHS}
=\langle Q^{-1/2}u,Q^{-1/2}v \rangle_{\SpHS}$. Then
$\mathcal{B}_{\SpQHS}=\{\hat{\xi}=Q^{1/2}\xi:\xi\in\mathcal{B}_{\SpHS}\}$
is an orthonormal basis of the Hilbert space $\SpQHS$.
\end{notation}

\begin{notation}\label{L20Space}
Let $L_2^0$ be the space of Hilbert-Schmidt operators from $\SpQHS$
to $\SpHS$ with the norm
\[
|\Phi|_{L_2^0}^2 =\sum_{\hat{\xi}\in\mathcal{B}_{{\SpQHS}}}
|\Phi\hat{\xi}|_{\SpHS}^2
=\sum_{\xi,\zeta\in\mathcal{B}_{\SpHS}}Q(\xi)|\langle \Phi\xi,\zeta
\rangle_{\SpHS}|^2 =\Tr [\Phi Q \Phi^\ast],
\]
where $Q(\xi)$ means $Q$ evaluated at $\xi$ as a positive function
on $\mathcal{B}_{\SpHS}$.
\end{notation}

\begin{lemma}\label{L20Norm}
If $\Psi\in L({\SpHS},{\SpHS})$, a bounded linear operator from
$\SpHS$ to $\SpHS$, then $\Psi$ restricted on $\SpQHS$ is a
Hilbert-Schmidt operator from $\SpQHS$ to $\SpHS$, and
$|\Psi|_{L_2^0}\leqslant \Tr(Q)\|\Psi\|^2$, where $\|\Psi\|$ is the
operator norm of $\Psi$.
\end{lemma}

\begin{proof}
\begin{align*}
|\Psi|_{L_2^0}^2 &=\sum_{\hat{\xi}\in\mathcal{B}_{\SpQHS}}
|\Psi\hat{\xi}|_{\SpHS}^2
\leqslant \|\Psi\|^2 \sum_{\hat{\xi}\in\mathcal{B}_{\SpQHS}}|\hat{\xi}|_{\SpHS}^2\\
&=\|\Psi\|^2 \sum_{\xi\in\mathcal{B}_{\SpHS}}\langle
Q^{1/2}\xi,Q^{1/2}\xi\rangle_{\SpHS} =\|\Psi\|^2
\sum_{\xi\in\mathcal{B}_{\SpHS}}\langle Q\xi,\xi\rangle_{\SpHS}
=\|\Psi\|^2 \Tr(Q)
\end{align*}
\end{proof}

\begin{notation}\label{BF}
Define $B:{\SpHS}\to L_2^0$ by $B(Y)A=(I+Y)A$ for $A\in \SpQHS$, and
$F:{\SpHS}\to {\SpHS}$ by $F(Y)=\frac{1}{2}(I+Y)D$.
\end{notation}
Note that $B$ is well--defined by Lemma \ref{L20Norm}. Also $D\in
{\SpHS}$, and so $F(Y)\in {\SpHS}$ and $F$ is well--defined as well.
\begin{theorem}\label{Main1}
The stochastic differential equation
\begin{align}
&dY_t = B(Y_t)dW_{t} + F(Y_t)dt \label{e.6.1}\\
&Y_0=0 \notag
\end{align}
has a unique solution, up to equivalence,  among the processes satisfying
\[
P\left(\int_0^T |Y_s|_{\SpHS}^2 ds < \infty \right) = 1.
\]
\end{theorem}

\begin{proof}
To prove this theorem we will use Theorem 7.4 from the book by
G.~DaPrato and J.~Zabczyk \cite{DaPratoBook1992} as it has been done
in \cite{Gordina2000a, Gordina2005a}. It is enough to check
\begin{enumerate}
\item[1.]
$B$ is a measurable mapping.
\item[2.]
$|B(Y_1)-B(Y_2)|_{L_2^0} \leqslant C_1 |Y_1-Y_2|_{\SpHS}$ for
$Y_1,Y_2\in {\SpHS}$;
\item[3.]
$|B(Y)|_{L_2^0}^2 \leqslant K_1(1+|Y|_{\SpHS}^2)$ for any $Y\in
{\SpHS}$;
\item[4.]
$F$ is a measurable mapping.
\item[5.]
$|F(Y_1)-F(Y_2)|_{\SpHS} \leqslant C_2 |Y_1-Y_2|_{\SpHS}$ for
$Y_1,Y_2\in {\SpHS}$;
\item[6.]
$|F(Y)|_{\SpHS}^2 \leqslant K_2(1+|Y|_{\SpHS}^2)$ for any $Y\in
{\SpHS}$.
\end{enumerate}
Proof of 1. By the proof of 2, $B$ is a continuous mapping,
therefore it is measurable.

\noindent
Proof of 2.
\begin{align*}
&|B(Y_1)-B(Y_2)|_{L_2^0}^2 =\sum_{\hat{\xi}\in\mathcal{B}_{\SpQHS}}
|(Y_1-Y_2)\hat{\xi}|_{\SpHS}^2
=\sum_{\xi\in\mathcal{B}_{\SpHS}} Q(\xi)|(Y_1-Y_2)\xi|_{\SpHS}^2\\
&\leqslant \sum_{\xi\in\mathcal{B}_{\SpHS}}
Q(\xi)\|\xi\|^2|Y_1-Y_2|_{\SpHS}^2 \leqslant
\max_{\xi\in\mathcal{B}_{\SpHS}}\|\xi\|^2
\left(\sum_{\xi\in\mathcal{B}_{\SpHS}}Q(\xi)\right)|Y_1-Y_2|_{\SpHS}^2
\\
& =\Tr Q\left(\max_{\xi\in\mathcal{B}_{\SpHS}}\|\xi\|^2\right)
|Y_1-Y_2|_{\SpHS}^2 = C_1^2 |Y_1-Y_2|_{\SpHS}^2,
\end{align*}
where $\|\xi\|$ is the operator norm of $\xi$, which is uniformly
bounded for all $\xi\in\mathcal{B}_{\SpHS}$.

\noindent
Proof of 3.
\begin{align*}
|B(Y_1)|_{L_2^0}^2
&=\sum_{\hat{\xi}\in\mathcal{B}_{\SpQHS}}|(I+Y)\hat{\xi}|_{\SpHS}^2
=\sum_{\xi\in\mathcal{B}_{\SpHS}} Q(\xi)|(I+Y)\xi|_{\SpHS}^2\\
&\leqslant |(I+Y)\xi|_{\SpHS}^2 \sum_{\xi\in\mathcal{B}_{\SpHS}}
Q(\xi)\|\xi\|^2 \le(1+|Y|_{\SpHS}^2)\cdot K_1.
\end{align*}

\noindent Proof of 4. By the proof of 5, $F$ is a continuous
mapping, therefore it is measurable.

\noindent
Proof of 5.
\begin{align*}
|F(Y_1)-F(Y_2)|_{\SpHS}=|\frac{1}{2}(Y_1-Y_2)D|_{\SpHS}
\le\|\frac{1}{2}D\||Y_1-Y_2|_{\SpHS}
\end{align*}

\noindent
Proof of 6.
\[
|F(Y)|_{\SpHS}^2=|\frac{1}{2}(I+Y)D|_{\SpHS}^2
\le\|\frac{1}{2}D\|^2|I+Y|_{\SpHS}^2 \leqslant K_2(1+|Y|_{\SpHS}^2).
\]
\end{proof}

\begin{notation}\label{BFSharp}
Let $B^\#:\SpHS \to L_2^0$  be the operator  $B^\#(Y)A=A^\#(I+Y)$,
and $F^\#: \SpHS \to \SpHS$ be the operator
$F^\#(Y)=\frac{1}{2}D^\#(Y+I)$.
\end{notation}

\begin{proposition}\label{YSharp}
If $Y_t$ is the solution to the stochastic differential equation
\begin{align*}
& dX_t = B(X_t)dW_{t} + F(X_t)dt\\
&X_0=0,
\end{align*}
where $B$ and $F$ are defined in Notation \ref{BF}, then $Y_t^\#$ is
the solution to the stochastic differential equation
\begin{align}
& dX_t = B^\#(X_t)dW_{t} + F^\#(X_t)dt \label{e.6.2}\\
&X_0=0, \notag
\end{align}
where $B^\#$ and $F^\#$ are defined in Notation \ref{BFSharp}.
\end{proposition}

\begin{proof}
This follows directly from the property $(AB)^\#=B^\# A^\#$ for any
$A$ and $B$, which can be verified by using part (5) of Proposition
\ref{SomeFacts1}.
\end{proof}

\begin{lemma}\label{Trace}
Let $U$ and $H$ be real Hilbert spaces. Let $\Phi:U\to H$ be a
bounded linear map. Let $G:H\to H$ be a bounded linear map. Then
\[
\Tr_H(G\Phi\Phi^\ast)=\Tr_U(\Phi^\ast G\Phi)
\]
\end{lemma}

\begin{proof}
\begin{align*}
\Tr_H(G\Phi\Phi^\ast)
&=\sum_{i,j\in H; k\in U} G_{ij}\Phi_{jk}(\Phi^\ast)_{ki}
=\sum_{i,j\in H; k\in U} G_{ij}\Phi_{jk}\Phi_{ik}\\
\Tr_U(\Phi^\ast G\Phi)
&=\sum_{i,j\in H; k\in U}(\Phi^\ast)_{ki}G_{ij}\Phi_{jk}
=\sum_{i,j\in H; k\in U} G_{ij}\Phi_{jk}\Phi_{ik}.
\end{align*}
Therefore $\Tr_H(G\Phi\Phi^\ast)=\Tr_U(\Phi^\ast G\Phi)$.
\end{proof}

\begin{lemma}\label{Sum_xi}
\[
\sum_{\xi\in\mathcal{B}_{\SpHS}}\big(Q^{1/2}\xi\big)\big(Q^{1/2}\xi\big)^\#=-D
\]
\end{lemma}

\begin{proof}
If $\xi\in\mathcal{B}_{\SpHS}$, then
$\xi\in\mathfrak{sp}\left(\infty\right)$, so $\xi^\#=-\xi$. We will
use the fact that
\[
(e_{ij}^{Re}e_{kl}^{Re})_{pq}=\delta_{ip}\delta_{jk}\delta_{lq}
\]
where $e_{ij}^{Re}$ is the matrix with the $(i,j)$th entry being $1$
and all other entries being zero. Using this fact, we see
\begin{enumerate}
\item
for
$\xi=\frac{1}{2}\big(e_{mn}^{Re}-e_{nm}^{Re}+e_{-m,-n}^{Re}
-e_{-n,-m}^{Re}\big)$ with $\sgn(mn)>0$,
\[
\big(Q^{1/2}\xi\big)\big(Q^{1/2}\xi\big)^\# =-\frac{1}{4}Q_{mn}^{Re}
\left[-e_{mm}^{Re}-e_{nn}^{Re}-e_{-m,-m}^{Re}-e_{-n,-n}^{Re}\right]
\]
\item
for
$\xi=\frac{1}{2}\big(e_{mn}^{Im}+e_{nm}^{Im}-e_{-m,-n}^{Im}
-e_{-n,-m}^{Im}\big)$ with $\sgn(mn)>0$,
\[
\big(Q^{1/2}\xi\big)\big(Q^{1/2}\xi\big)^\# =-\frac{1}{4}Q_{mn}^{Im}
\left[-e_{mm}^{Re}-e_{nn}^{Re}-e_{-m,-m}^{Re}-e_{-n,-n}^{Re}\right]
\]
\item
for
$\xi=\frac{1}{2}\big(e_{mn}^{Re}+e_{-n,-m}^{Re}+e_{-m,-n}^{Re}
+e_{n,m}^{Re}\big)$ with $\sgn(mn)<0$,
\[
\big(Q^{1/2}\xi\big)\big(Q^{1/2}\xi\big)^\# =-\frac{1}{4}Q_{mn}^{Re}
\left[e_{mm}^{Re}+e_{nn}^{Re}+e_{-m,-m}^{Re}+e_{-n,-n}^{Re}\right]
\]
\item
for
$\xi=\frac{1}{2}\big(e_{mn}^{Im}+e_{-n,-m}^{Im}-e_{-m,-n}^{Im}
-e_{nm}^{Im}\big)$ with $\sgn(mn)<0$,
\[
\big(Q^{1/2}\xi\big)\big(Q^{1/2}\xi\big)^\# =-\frac{1}{4}Q_{mn}^{Im}
\left[e_{mm}^{Re}+e_{nn}^{Re}+e_{-m,-m}^{Re}+e_{-n,-n}^{Re}\right].
\]
\end{enumerate}

Each of the above is a diagonal matrix. The lemma can be proved by
looking at the diagonal entries of the sum.
\end{proof}

\begin{theorem}\label{Main2}
Let $Y_t$ be the solution to Equation \ref{e.6.1}. Then $Y_t+I\in
\Sp(\infty)$ for any $t>0$ with probability $1$.
\end{theorem}

\begin{proof}
The proof is adapted from papers by M.~Gordina \cite{Gordina2000a,
Gordina2005a}. Let $Y_t$ be the solution to Equation \eqref{e.6.1}
and $Y_t^\#$ be the solution to Equation \eqref{e.6.2}. Consider the
process $\textbf{Y}_t=(Y_t,Y_t^\#)$ in the product space
${\SpHS}\times {\SpHS}$. It satisfies the following stochastic
differential equation
\[
d\textbf{Y}_t=(B(Y_t),B^\#(Y_t^\#))dW + (F(Y_t),F^\#(Y_t^\#))dt.
\]

Let $G$ be a function on the Hilbert space ${\SpHS}\times {\SpHS}$
defined by $G(Y_1,Y_2)=\Lambda((Y_1+I)(Y_2+I))$, where $\Lambda$ is
a nonzero linear real bounded functional from ${\SpHS}\times
{\SpHS}$ to $\mathbb{R}$. We will apply It\^o's formula to
$G(\textbf{Y}_t)=G(Y_t,Y_t^\#)$. Then $(Y_t+I)(Y_t^\#+I)=I$ if and
only if $\Lambda((Y_t+I)(Y_t^\#+I)-I)=0$ for any $\Lambda$.

In order to use It\^o's formula we must verify that $G$ and the
derivatives $G_t$, $G_{\textbf{Y}}$, $G_{\textbf{YY}}$ are uniformly
continuous on bounded subsets of $[0,T]\times {\SpHS}\times
{\SpHS}$, where $G_{\textbf{Y}}$ is defined as follows
\[
G_{\textbf{Y}}(\textbf{Y})(\textbf{S}) =\lim_{\epsilon\to0}
\frac{G(\textbf{Y}+\epsilon\textbf{S})-G(\textbf{Y})}{\epsilon}
\hspace{.2in}\mbox{for any }\textbf{Y},\textbf{S}\in {\SpHS}\times
{\SpHS}
\]
and $G_{\textbf{YY}}$ is defined as follows
\[
G_{\textbf{YY}}(\textbf{Y})(\textbf{S}\otimes\textbf{T})
=\lim_{\epsilon\to0}
\frac{G_{\textbf{Y}}(\textbf{Y}+\epsilon\textbf{T})(\textbf{S})
-G_{\textbf{Y}}(\textbf{Y})(\textbf{S})}{\epsilon}
\]
for any $\textbf{Y},\textbf{S},\textbf{T}\in {\SpHS}\times {\SpHS}$.
Let us calculate $G_t$, $G_{\textbf{Y}}$, $G_{\textbf{YY}}$.
Clearly, $G_t=0$. It is easy to verify that for any
$\textbf{S}=(S_1,S_2)\in {\SpHS}\times {\SpHS}$
\[
G_{\textbf{Y}}(\textbf{Y})(\textbf{S})
=\Lambda(S_1(Y_2+I)+(Y_1+I)S_2)
\]
and for any $\textbf{S}=(S_1,S_2)\in {\SpHS}\times {\SpHS}$ and
$\textbf{T}=(T_1,T_2)\in {\SpHS}\times {\SpHS}$
\[
G_{\textbf{YY}}(\textbf{Y})(\textbf{S}\otimes\textbf{T})
=\Lambda(S_1T_2+T_1S_2).
\]
So the condition is satisfied.

We will use the following notation
\begin{align*}
&G_{\textbf{Y}}(\textbf{Y})(\textbf{S})
=\langle \bar{G}_{\textbf{Y}}(\textbf{Y}),\textbf{S} \rangle_{\SpHS\times\SpHS}\\
\noalign{\vskip .05 true in}
&G_{\textbf{YY}}(\textbf{Y})(\textbf{S}\otimes\textbf{T})= \langle
\bar{G}_{\textbf{YY}}(\textbf{Y})\textbf{S},\textbf{T}
\rangle_{\SpHS\times\SpHS},
\end{align*}
where $\bar{G}_{\textbf{Y}}(\textbf{Y})$ is an element of
${\SpHS}\times {\SpHS}$ corresponding to the functional
$G_{\textbf{Y}}(\textbf{Y})$ in $({\SpHS}\times {\SpHS})^\ast$ and
$\bar{G}_{\textbf{YY}}(\textbf{Y})$ is an operator on ${\SpHS}\times
{\SpHS}$ corresponding to the functional
$G_{\textbf{YY}}(\textbf{Y})\in (({\SpHS}\times
{\SpHS})\otimes({\SpHS}\times {\SpHS}))^\ast$.

Now we can apply It\^o's formula to $G(\textbf{Y}_t)$
\begin{align*}
G(\textbf{Y}_t)-G(\textbf{Y}_0)
=&\int_0^t \langle \bar{G}_{\textbf{Y}}(\textbf{Y}_s),
\big(B(Y_s)dW_s,B^\#(Y_s^\#)dW_s\big)\rangle_{\SpHS\times\SpHS}\\
+&\int_0^t \langle \bar{G}_{\textbf{Y}}(\textbf{Y}_s),
\big(F(Y_s),F^\#(Y_s^\#)\big)\rangle_{\SpHS\times\SpHS} ds\\
+&\int_0^t\frac{1}{2}\Tr_{{\SpHS}\times
{\SpHS}}\bigg[\bar{G}_{\textbf{YY}}(\textbf{Y}_s)
\Big(B(Y_s)Q^{1/2},B^\#(Y_s^\#)Q^{1/2}\Big)\\
&\hspace{1.5in}\Big(B(Y_s)Q^{1/2},B^\#(Y_s^\#)Q^{1/2}\Big)^\ast\bigg]
ds.
\end{align*}

Let us calculate the three integrands separately.
The first integrand is
\begin{align*}
\langle \bar{G}_{\textbf{Y}}(\textbf{Y}_s),
&\big(B(Y_s)dW_s,B^\#(Y_s^\#)dW_s\big)\rangle_{\SpHS\times\SpHS}\\
&=\Big(B(Y_s)dW_s\Big)(Y_s^\#+I)+(Y_s+I)\Big(B^\#(Y_s^\#)dW_s\Big)\\
&=(Y_s+I)dW_s(Y_s^\#+I)+(Y_s+I)dW_s^\#(Y_s^\#+I)
=0.
\end{align*}

The second integrand is
\begin{align*}
\langle \bar{G}_{\textbf{Y}}(\textbf{Y}_s),
&\big(F(Y_s),F^\#(Y_s^\#)\big)\rangle_{\SpHS\times\SpHS}\\
&=F(Y_s)(Y_s^\#+I)+(Y_s+I)F^\#(Y_s^\#)\\
&=\frac{1}{2}(Y_s+I)D(Y_s^\#+I)+\frac{1}{2}(Y_s+I)D^\#(Y_s^\#+I)\\
&=\frac{1}{2}(Y_s+I)(D+D^\#)(Y_s^\#+I)\\
&=(Y_s+I)D(Y_s^\#+I),
\end{align*}
where we have used the fact that $D=D^\#$, since $D$ is a diagonal matrix
with all real entries.

The third integrand is
\begin{align*}
& \frac{1}{2}\Tr_{{\SpHS}\times {\SpHS}}
\\
& \left[\bar{G}_{\textbf{YY}}(\textbf{Y}_s)
\left(B(Y_s)Q^{1/2},B^\#(Y_s^\#)Q^{1/2}\right)
\left(B(Y_s)Q^{1/2},B^\#(Y_s^\#)Q^{1/2}\right)^\ast\right]
\\
&
=\frac{1}{2}\Tr_{{\SpHS}}\left[\left(B(Y_s)Q^{1/2},B^\#(Y_s^\#)Q^{1/2}\right)^\ast
\bar{G}_{\textbf{YY}}(\textbf{Y}_s)
\left(B(Y_s)Q^{1/2},B^\#(Y_s^\#)Q^{1/2}\right) \right]\\
&=\frac{1}{2}\sum_{\xi\in\mathcal{B}_{\SpHS}}G_{\textbf{YY}}(\textbf{Y}_s)
\bigg(\Big(B(Y_s)Q^{1/2}\xi,B^\#(Y_s^\#)Q^{1/2}\xi\Big)\\
&\hspace{2in}
\otimes\Big(B(Y_s)Q^{1/2}\xi,B^\#(Y_s^\#)Q^{1/2}\xi\Big)\bigg)\\
&=\sum_{\xi\in\mathcal{B}_{\SpHS}}
\Big(B(Y_s)Q^{1/2}\xi\Big)\Big(B^\#(Y_s^\#)Q^{1/2}\xi\Big)\\
&=\sum_{\xi\in\mathcal{B}_{\SpHS}}(Y_s+I)
\bigg(\big(Q^{1/2}\xi\big)\big(Q^{1/2}\xi\big)^\#\bigg)(Y_s^\#+I)\\
&=-(Y_s+I)D(Y_s^\#+I),
\end{align*}
where the second equality follows from Lemma \ref{Trace}, and the
last equality follows from Lemma \ref{Sum_xi}.

The above calculations show that the stochastic differential of $G$
is zero. So $G(\textbf{Y}_t)=G(\textbf{Y}_0)=\Lambda(I)$ for any
$t>0$ and any nonzero linear real bounded functional $\Lambda$ on
${\SpHS}\times {\SpHS}$. This means $(Y_t+I)(Y_t^\#+I)=I$ almost
surely for any $t>0$. Similarly we can show $(Y_t^\#+I)(Y_t+I)=I$
almost surely for any $t>0$. Therefore $Y_t+I\in \Sp(\infty)$ almost
surely for any $t>0$.
\end{proof}

\bibliographystyle{amsplain}

\end{document}